\documentclass[12pt, reqno]{amsart}
\usepackage{amscd,amsthm,amsfonts,amssymb,amsmath,euscript}
\usepackage[dvips]{graphicx}
\usepackage[dvips]{graphics}
\usepackage[matrix,arrow]{xy}



\newcounter{statements}

\newtheorem{theorem}[statements]{Theorem}
\newtheorem{conjecture}[statements]{Conjecture}
\newtheorem{problem}[statements]{Problem}
\newtheorem{question}[statements]{Question}

\newtheorem{proposition}[statements]{Proposition}

\theoremstyle{definition}

\newtheorem{definition}[statements]{Definition}

\theoremstyle{remark}
\newtheorem{remark}[statements]{Remark}

\newcommand{\QQ}{{\mathbb Q}}
\newcommand{\ZZ}{{\mathbb Z}}
\newcommand{\NN}{{\mathbb N}}

\newcommand{\PP}{{\mathbb P}}
\newcommand{\CC}{{\mathbb C}}

\newcommand{\Proof}{{\bf Proof}}

\newfont{\smallskob}{cmbx7 scaled\magstep4}
\newfont{\bigskob}{cmbx12 scaled\magstep4}

\newcommand{\pic}{\mathrm{Pic}\,}

\newcommand{\wt}{\mathrm{wt}}
\newcommand{\cO}{\mathcal{O}}

\newcommand{\tit}{Hori--Vafa mirror models for complete
intersections in weighted projective spaces and weak
Landau--Ginzburg models}

\begin{document}

\begin{title}
\tit
\end{title}

\author{Victor Przyjalkowski}

\thanks{The work was partially supported by FWF grant P20778,
RFFI grants 11-01-00336-a and 11-01-00185-a, grants
NSh$-4713.2010.1$,  MK-503.2010.1, and AG Laboratory GU-HSE, RF government
grant, ag. 11 11.G34.31.0023.}

\address{Steklov Mathematical Institute, 8 Gubkina street, Moscow 119991, Russia} %

\email{victorprz@mi.ras.ru, victorprz@gmail.com}

\maketitle

\begin{abstract}
We prove that Hori--Vafa mirror models for smooth Fano complete
intersections in weighted projective spaces admit an interpretation
as Laurent polynomials.
\end{abstract}





Mirror symmetry
of variations of Hodge structures 
 states that for any smooth
Fano variety $X$ there exists a dual Landau--Ginzburg model $f\colon
Y\to \CC$ such that an essential part of the regularized quantum
differential equation for $X$ is of Picard--Fuchs type. In other
words, the solutions of a certain differential equation (constructed
via genus 0 Gromov--Witten invariants for $X$ --- the numbers which
count rational curves lying on $X$) are the periods of the dual
family
(for more details and references see~\cite{Prz09a}). By
definition, the
relevant Picard--Fuchs 
differential equation depends only on relative birational type of
$Y$. If one assumes $Y=(\CC^*)^N$ one can translate mirror
correspondence to the quantitative level, that is, to combinatorics
of Laurent polynomials. Then $f$ may be represented by a Laurent
polynomial, which is called a (very) weak Landau--Ginzburg model.
The following conjecture states that this hypothesis is not
very restrictive, particulary for the 
case of $\pic X=\ZZ$.

\begin{conjecture}[\cite{Prz09a}]
\label{conjecture} Any smooth Fano variety of dimension $N$ with
Picard rank $1$ has a weak Landau--Ginzburg model $f\in \CC[x_1^{\pm
1},\ldots, x_N^{\pm 1}]$.
\end{conjecture}

This conjecture holds for threefolds, complete intersections in projective spaces,
Grassmannians and some complete intersections therein, varieties
admitting degenerations to``good'' toric varieties (for more details see~\cite{Prz09a}). In the paper we
prove that it also holds for smooth complete intersections of
Cartier divisors in weighted projective spaces. That is, we prove
that Hori--Vafa suggestions for Landau--Ginzburg models for such
varieties may be interpreted as Laurent polynomials.
The similar problems (for weighted projective spaces in an orbifold setup) are studied in~\cite{Dou06} and~\cite{DM09}.

\section{Hori-Vafa models}
We give some basic definitions and notions about weighted projective spaces
and complete intersections therein mostly following~\cite{Do82}.

We consider weighted projective spaces as projective varieties (not as smooth stacks). We denote by $(a_1,\ldots,a_r)$ the greatest common divisor of
$a_1,\ldots,a_r\in \NN$.

\begin{definition}
A weighted projective space $\PP(w_0, \ldots, w_n)$  is called
\emph{normalized} if for any $i$ we have
$(w_0,\ldots,\hat{w}_i,\ldots,w_n)=1$ and $w_0\leq w_1\leq
\ldots\leq w_n$.
\end{definition}

\begin{remark}
It is easy to see that $\PP(w_0, dw_1 \ldots, dw_n)\cong \PP(w_0,w_1
\ldots, w_n)$, so any weighted projective space is isomorphic to a
unique normalized one.
\end{remark}

\begin{definition}
The zero set of (weighted) homogenous polynomial $f\in
\CC[x_0,\ldots,x_n]$, $\wt (x_i)=w_i$, of weighted degree $d$ is
called \emph{a hypersurface of degree $d$} in
$\PP=\PP(w_0,\ldots,w_n)$.
\end{definition}

As the rank of the Weil group of a weighted projective space is 1,
any effective 
Weil divisor is proportional to the zero locus of some weighted
homogenous polynomial. Its degree is called the degree of the
divisor. It is easy to see that a Weil divisor of degree $d$ is
Cartier if and only if $d$ is integral and all $w_i$'s divide $d$.

The singular locus of $\PP$ is the union of subvarieties of form
$\{x_{i_1}=\ldots=x_{i_r}=0\}$, where $w_{i_1},\ldots,w_{i_r}$ is a
minimal collection of weights such that the rest of the weights have
common
prime divisor. 
Consider a complete intersection $X=X_1\cap
\ldots \cap X_k$, where $X_1,\ldots,X_k$ are Cartier divisors.
It is quasismooth as a complete intersection of weighted Fermat hypersurfaces is quasismooth.
By Proposition~8 in~\cite{Di86} together with Proposition~2 in loc.\,cit.
the singularities of $X$ are
the intersection of $X$ with the singularities of $\PP$.
In particular $X$ is
smooth if and only if the maximal dimension of the strata of
singularities of $\PP$ is less then $k$. This means that
$(w_{i_1},\ldots,w_{i_{k+1}})=1$ for any collection of weights
$w_{i_1},\ldots,w_{i_{k+1}}$ (cf.~\cite{Di86}).

Let $\deg X_i=d_i$. The canonical sheaf of $X$ is
$\cO(d_1+\ldots+d_k-w_0-\ldots-w_n)|_X$. So $X$ is Fano if and only
if $\sum d_i<\sum w_j$.

\begin{definition}[\cite{HV00}, see also~\cite{Giv96}]
Consider a smooth complete intersection $X=X_1\cap\ldots\cap X_k$ in
$\PP(w_0,\ldots,w_n)$ such that $X_i$ is a Cartier divisor of degree
$d_i$ and there are $k$ non-intersecting subsets $I_i\subset
\{0,\ldots,n\}$ such that $\sum_{j\in I_i} {w_j}=d_i$ (we call this splitting a
\emph{$\QQ$-nef-partition\footnote{It is called \emph{a nef-partition} in Gorenstein case.}}
). Then \emph{a Hori--Vafa model} for $X$ is an affine
variety
$$
\left\{%
\begin{array}{l}
    x_0^{w_0}\cdot\ldots \cdot x_n^{w_n}=1, \\
    \sum_{j\in I_i} {x_j}=1 \\
\end{array}%
\right.
$$
with function (potential) $f=x_0+\ldots+x_n$. Up to a shift
$f\mapsto f-k$ we can define the potential as the sum of variables
whose indices do not lie in any $I_i$'s.
\end{definition}

\section{Weak Landau--Ginzburg models}

We define a (very) weak Landau--Ginzburg models for Fano variety
following~\cite{Prz08}.

Let $X$ be a smooth Fano variety. Given its Gromov--Witten
invariants (the numbers which count rational curves lying on $X$)
one can construct the so called~\emph{regularized quantum
differential equation} for $X$. This equation for a complete intersection in weighted projective spaces is \emph{of type DN} and has a unique
normalized analytic solution $I^X_{H^0}(t)$ called \emph{the
constant term of regularized $I$-series}.

Consider a Laurent polynomial $f\in \CC[x_1^{\pm 1}, \ldots,x_n^{\pm
1}]$. Let $b_i$ be the constant term of $f^i$. The series
$\Phi_f(t)=1+b_1t+b_2t^2+\ldots$ is called \emph{the constant terms
series} for $f$. It is 
an analytic solution of
the Picard--Fuchs differential equation for a pencil of
hypersurfaces in the torus given by $f$.

\begin{definition}
The polynomial $f$ (or the pencil of hypersurfaces associated to it)
is called \emph{a very weak Landau--Ginzburg model} for $X$ if the
regularized quantum differential equation for $X$ coincides with the
Picard--Fuchs equation for $f$, or, equivalently, if
$I^X_{H^0}=\Phi_f$.

It is called \emph{a weak Landau--Ginzburg model} for $X$ if, in
addition, the general element of the pencil given by $f$ is
birational to a Calabi--Yau variety.
\end{definition}

\begin{proposition}
Let $X$ be a smooth complete intersection of Cartier divisors of
degrees $d_1,\ldots,d_k$ in normalized weighted projective space
$\PP(w_0,\ldots,w_n)$. Assume that $X$ is a Fano variety. Then there
are $k$ non-intersecting subsets $I_1, \ldots,I_k\subset
\{0,\ldots,n\}$ such that $d_i=\sum_{j\in I_i} w_j$ for any $i$ and
$w_j=1$ for all $j\notin I_1\cup\ldots\cup I_k$. \label{proposition}
\end{proposition}

\Proof. 
We have the following numerical conditions:
$$
w_0\leq w_1\leq \ldots\leq w_n,
$$
$$
w_j|d_i,\ \ \ \ j=0\ldots n,\ i=1,\ldots, k,
$$
$$
\sum d_i<\sum w_i,
$$
$$
(w_{i_1},\ldots,w_{i_{k+1}})=1,\ \ \ \ \{i_1,\ldots,i_{k+1}\}\subset
\{0,\ldots,n\}.
$$
%
%
Apply the following ``reduction process''.
Let $p$ be a divisor of one of the weights. Divide by $p$ all the
degrees and those of the weights which are divisible by $p$. Up
to renumbering of the weights we get a collection of weights and
degrees satisfying the conditions above. Repeat the procedure until
all weights become equal to 1. Consider $k$ non-intersecting subsets
$I_1,\ldots,I_k$ of $\{0,\ldots,n\}$ whose orders equal the degrees
we got on the last step. Let $I=\{0,\ldots,n\}\setminus
\{I_1\cup\ldots\cup I_k\}$.

Denote $\sum w_j$, $w_j\in I_i$, by $|I_i|$. Start the reduction
process in the reverse direction. The weights and the degrees change
on each step. Change the elements of $I$ and $I_i$'s on the first
step in such a way that each of $I_i$'s contains at most one index
of an increasing weight and $I$ contains none of them. Change $I$
and $I_i$'s on each step in the following way: if $|I_i|=d_i$ (where
$d_i$'s are the degrees on this step) do nothing. If $|I_i|<d_i$,
add indices from $I$ whose weights increase to $I_i$ (it is easy to
check that $d_i-|I_i|$ is not less then the prime divisor $p$ we
increase). As the number of increasing weights is not greater than
$k$ we get $I$ containing only indices corresponding to the weight
1. We get $|I_i|\leq d_i$. If $|I_i|<d_i$, add $d_i-|I_i|$ indices
from $I$. Doing such changes of $I$ and $I_i$'s on each step we get
$|I_i|=d_i$, and all the weights whose indices lie in $I$ equal 1.
Finally we get the partition we need. \qed

\medskip

\begin{remark}
Let $d_0=\sum w_i-\sum d_j$ be the Fano index of $X$. It is easy to
see from the proof of Proposition~\ref{proposition} that there are
actually at least $d_0+1$ weights that are equal to 1. This bound is
strict. The example is hypersurface of degree 6 in $\PP(1,1,2,3)$.
\end{remark}

\begin{theorem}
Let $X$ be a smooth Fano complete intersection of Cartier divisors
in normalized weighted projective space $\PP(w_0,\ldots,w_n)$. Then
$X$ has a very weak Landau--Ginzburg model. \label{theorem}
\end{theorem}

\Proof. In order to keep the notation simple we prove this theorem
for the case when $X$ is a hypersurface; the general case can be proved
identically. Prove that the pencil given by the Hori--Vafa model for
$X$ is relative birationally isomorphic to a pencil of hypersurfaces
in $(\CC^*)^{n-1}$. By Proposition~\ref{proposition}, $d_1$ is the
sum of some weights such that the rest of the weights equal 1.
Renumber weights for convenience such that $d_1=w_{r+1}+\ldots+w_n$.
Then do the well-known trick with a projective change of coordinates
for Hori--Vafa model for $X$. That is, the Hori--Vafa model is
$$
\left\{%
\begin{array}{l}
    x_0^{w_0}\cdot\ldots \cdot x_n^{w_n}=1, \\
    x_{r+1}+\ldots+x_n=1 \\
\end{array}%
\right.
$$
with potential $f=x_0+\ldots+x_r$. Consider (projective) change of
coordinates
$$
x_i=\frac{y_i}{y_{r+1}+\ldots+y_n}, \ \ \ \ i=r+1,\ldots,n.
$$
The second equation of the system disappears. As $w_0=1$ we may
rewrite the first variable as
$$
x_0=\frac{(y_{r+1}+\ldots+y_n)^{d_1}}{x_1^{w_1}\cdot\ldots\cdot
x_r^{w_r} y_{r+1}^{w_{r+1}}\cdot\ldots\cdot y_n^{w_n}}.
$$
In the local chart, say, $y_n=1$ we finally get a very weak
Landau--Ginzburg model
$$
f_X=\frac{(y_{r+1}+\ldots+y_{n-1}+1)^{d_1}}{x_1^{w_1}\cdot\ldots\cdot
x_r^{w_r} \cdot y_{r+1}^{w_{r+1}}\cdot\ldots\cdot
y_{n-1}^{w_{n-1}}}+x_1+\ldots+x_r.
$$

The constant term of regularized $I$-series for $X$ is given by
$$
I^X_{H^0}=\sum_{m=0}^{\infty}
\frac{(d_0m)!(d_1m)!}{(w_0m)!\cdot\ldots\cdot (w_nm)!}t^{d_0m}
$$
(see~\cite{Prz07c}). One can see that this series coincides with the
constant terms series for $f_X$. \qed

%
%
%


\begin{remark}
Let $X_d$ be a smooth hypersurface in a weighted projective space
such that its Fano index $i_X=w_0+\ldots+w_n-d$ is one. Let $f_X$ be
a very weak Landau--Ginzburg model $f_X$ given by
Theorem~\ref{theorem}. One can prove (see the proof of Theorem 14
in~\cite{Prz09a}) that a general element of the pencil defined by
$f_X$ birational to a Calabi--Yau variety. In other words, $f_X$ is
actually a weak Landau--Ginzburg model for $X$.
\end{remark}

\begin{problem}
Prove that this is true for all smooth Fano complete intersections
of Cartier divisors in weighted projective spaces. \label{problem1}
\end{problem}

%
%
%

\begin{remark}
\label{toric}
In~\cite{IP11} N.\,Ilten and the author prove that very
weak Landau--Ginzburg models of Hori--Vafa type for complete intersections are toric. This means that their Newton polytopes are fan polytopes
of toric degenerations of these complete intersections.
\end{remark}

It seems that the assumptions on varieties we need for Hori--Vafa
mirror models can be weakened. One can consider a complete
intersection $X$ which does not intersect the singular locus of
$\PP$ (in a classical setup this condition is necessary, because otherwise $X$ should be
considered as an orbifold). 
On the numerical level this means that for any $q>1$ the number of
weights divisible by $q$ is not greater then the number of degrees
divisible by $q$. In all examples we consider we still get an
appropriate $\QQ$-nef-partition. 

\begin{question}
Is it always true? If not, what conditions should we put to have an appropriate $\QQ$-nef partition?\label{question:quasi-projective}
\end{question}

\begin{remark}
If this is true then this statement can be strengthened in the following way. There is a $\QQ$-nef-partition $I_1,\ldots,I_k$ such that
$w_i=1$ for $i\notin I_1\cap\ldots\cap I_k$. \label{collary}
%
Indeed, consider the given $\QQ$-nef-partition. 
Delete $w_0$. 
All numerical conditions still hold for a collection of weights and degrees we get.
Thus there is another
appropriate $\QQ$-nef-partition. Deleting the smallest weights step by
step we
find the partition we need. 
\end{remark}


It is also natural to consider Hori--Vafa models for quasismooth
Fano complete intersections. But even in the case of a Cartier
hypersurface it is not always possible to write down a Hori--Vafa model.
An example is a hypersurface of degree 30 in $\PP(1, 6, 10, 15)$: it has no $\QQ$-nef partition.
Another example, suggested to the author by S.\,Galkin, is a hypersurface of degree 30 in $P(1, 6, 6, 6, 6, 10, 10, 15)$.
It shows that even nef-partition (that is $\mathbb Q$-nef-partition in Gorenstein variety) does not necessarily exist.

The reason of this phenomenon seems to be the following: such complete intersections
should be considered as \emph{smooth stacks} instead of considering them as \emph{singular varieties}.

\begin{question}
Is there a stacky version of Hori--Vafa procedure? If yes, can it be reformulated in Laurent polynomials terms?
\end{question}

%

Even if a hypersurface had a Hori--Vafa model, it can
have no very weak Landau--Ginzburg model of type discussed in the paper. An example is a
hypersurface of degree 30 in $\PP(1,1,1,1,1,6,10,15)$.

\begin{question}
Does this hypersurface (or all complete intersections having Hori--Vafa models) have another weak Landau--Ginzburg models, not of Hori--Vafa type?
In other words, is it rational?
\end{question}

\medskip

The author is grateful to I.\,Cheltsov, S.\,Galkin, V.\,Golyshev, L.\,Katzarkov, D.\,Orlov,
K.\,Shramov, D.\,Stepanov, and D.\,van Straten for helpful comments and
T.\,Logvinenko for detailed proofreading and English checking.


\end{document}